\documentclass[12pt]{article}
\usepackage[top=4cm,bottom=3cm,left=4cm,right=3cm]{geometry}
\usepackage[round]{natbib}

\usepackage{amssymb}
\usepackage{amsthm}

\usepackage{tabularx,environ}
\usepackage{soul}
\usepackage{caption}
\usepackage{subcaption}
\usepackage{tikz}
\usepackage{amsmath}
\usepackage{multirow}
\usepackage{color}

\newtheorem{proposition}{Proposition}

\usepackage{rotating}
\usepackage{fixmath}
\usepackage{algorithm}
\usepackage[noend]{algpseudocode}
\usepackage{lscape}
\usepackage{hyperref}
\usepackage{amssymb}
\usepackage[makeroom]{cancel}
\usepackage[affil-it]{authblk} 
\makeatletter
\newcommand{\problemtitle}[1]{\gdef\@problemtitle{#1}}% Store problem title
\newcommand{\probleminput}[1]{\gdef\@probleminput{#1}}% Store problem input
\newcommand{\problemquestion}[1]{\gdef\@problemquestion{#1}}% Store problem question

\NewEnviron{proble}{
  \problemtitle{}\probleminput{}\problemquestion{}% Default input is empty
  \BODY% Parse input
  \par\addvspace{.5\baselineskip}
  \noindent
  \begin{tabularx}{\textwidth}{@{\hspace{\parindent}} l X c}
    \multicolumn{2}{@{\hspace{\parindent}}l}{\@problemtitle} \\% Title
    \textbf{Input:} & \@probleminput \\% Input
    \textbf{Question:} & \@problemquestion% Question
  \end{tabularx}
  \par\addvspace{.5\baselineskip}
}
\newcommand\notsotiny{\@setfontsize\notsotiny{6.3}{6.3}}
\makeatother

\begin{document}

%\begin{frontmatter}

\title{Recovering feasibility in real-time conflict-free vehicle routing}

% \author[1]{XXX\thanks{corresponding author}}
\author[1]{Tommaso Adamo}
%\ead{tommaso.adamo@unisalento.it}

\author[1]{Gianpaolo Ghiani\thanks{corresponding author}}%\corref{mycorrespondingauthor}}
%\ead{gianpaolo.ghiani@unisalento.it}

\author[1]{Emanuela Guerriero}
%\ead{emanuela.guerriero@unisalento.it}

%\address[salento]{Dipartimento di Ingegneria dell'Innovazione,  Universit\`{a} del Salento, \\Via per Monteroni, 73100 Lecce, Italy}
\affil[1]{Dipartimento di Ingegneria dell'Innovazione,  Universit\`{a} del Salento }

\maketitle
\begin{abstract}
\textit{Conflict-Free Vehicle Routing Problems} (CF-VRPs) arise in manufacturing, transportation and logistics facilities where \textit{Automated Guided Vehicles} (AGVs) are utilized to move loads. Unlike \textit{Vehicle Routing Problems} arising in distribution management, CF-VRPs explicitly consider the limited capacity of the arcs of the guide path network to avoid collisions among vehicles. AGV applications have two peculiar features. First, the uncertainty affecting both travel times and machine ready times may result in vehicle delays or anticipations with respect to the fleet nominal plan. Second, the relatively high vehicle speed (in the order of one or two meters per second) requires vehicle plans to be revised in a very short amount of time (usually few milliseconds) in order to avoid collisions. In this paper we present fast exact algorithms to recover plan feasibility in real-time. In particular, we identify two corrective actions that can be implemented in real-time and formulate the problem as a linear program with the aim to optimize four common performance measures (total vehicle delay, total weighted delay, maximum route duration and total lateness). Moreover, we develop tailored algorithms which, tested on randomly generated instances of various sizes, prove to be three orders of magnitude faster than using off-the-shelf solvers.
\end{abstract}


%\begin{keyword}
%Logistics \sep Automated Guided Vehicles \sep Conflict-free Vehicle Routing and Scheduling \sep Real-time Optimization
%\end{keyword}
\providecommand{\keywords}[1]{\textbf{\textit{Logistics \sep Automated Guided Vehicles \sep Conflict-free Vehicle Routing and Scheduling \sep Real-time Optimization}} #1}
%\end{frontmatter}

\section{Introduction} \label{intro}
\textit{Conflict-Free Vehicle Routing and Scheduling Problems} (CF-VRPs) arise in manufacturing plants \citep{ullrich2015automated}, warehouses \citep{van1999models} and automated port terminals \citep{stahlbock2008vehicle, schwientek2017literature} where driverless vehicles (usually referred to as Automated Guided Vehicles, AGVs) are used to move materials, components and finished products, often in the form of palletized and containerized loads. 
Independently of the guidance system (e.g., magnet spot, laser or GPS navigation), an AGV may be assumed to move along the arcs of a graph whose vertex set represents loading/unloading stations, storage positions as well as intersections of segments of the guide path network. A distinctive feature of these systems is that, in order to avoid collisions among vehicles, a number of non-overlapping constraints have to be imposed. For instance, one may require that at most one vehicle occupies an arc at any given time \textit{or} that vehicles must keep a minimum distance, etc.
Relevant contributions to the solution of CF-VRPs have been proposed by \cite{krishnamurthy1993developing}, \cite{desaulniers2003dispatching}, \cite{correa2007scheduling} and \cite{miyamoto2016local}. More recently, \cite{adamo2018path} have studied the problem of determining vehicle routes (and vehicle speeds on routes segments) in such a way that no conflict arise, time windows are met, and the total energy consumption is minimized. 

In most industrial applications, arc travel times and machine ready times are not known exactly in advance since they depend, to some extent, on a number of variables that cannot be controlled (e.g., enable signals from automatic machines, battery levels, dirt on the floor, etc.) or on unpredictable events (e.g., a worker cutting the road to an AGV, a machine breakdown). Consequently, it may happen that an initially feasible CF-VRP solution (fleet \textit{nominal plan}) becomes unfeasible as a result of the delays (or anticipations) accumulated by some vehicles. The problem studied in this paper is how to recover plan feasibility after such variations have occurred. The objective is to avoid conflicts by implementing some \textit{corrective actions}. We consider two corrective actions, namely (a) imposing a delay to some vehicles and (b) accelerating some vehicles. The aim is to minimize four of the most common performance measures: (1) total vehicle delay, (2) total weighted delay, (3) maximum route duration (\textit{makespan}) and (4) total lateness w.r.t. specified \textit{due dates}.

Before delving into the description of our contribution, we position our paper towards a number of related research lines. The need to coordinate the routes and schedules of a fleet of vehicles to maintain a minimum distance from each other at any time arises in several sectors, including air traffic control and railway traffic management. In \textit{air traffic control}, aircraft move in a shared airspace and can not get closer to each other than a given safety distance in order to avoid possible conflicts. \cite{ribeiro2020review} review  conflict resolution methods for manned and unmanned aviation. In this line of research, a significant contribution is proposed by \cite{pallottino2002conflict} that consider the path planning problem among given waypoints avoiding all possible conflicts with the objective to minimize the total flight time. The authors propose two formulations as a mixed-integer linear program (MILP): in the first, only velocity changes are admissible maneuvers while in the second only heading angle changes are allowed. Solutions were obtained with standard optimization software. Unlike our contribution, this line of research assumes that aircraft move in a continuous (Euclidean) airspace. Another difference is that several seconds are usually allowed to generate a feasible solution.

In \textit{railway traffic management}, initial relatively small delays of some trains may propagate and generate conflicts. As a result, especially during congested traffic situations where the infrastructure capacity is completely exploited for trains circulation, some trains must be stopped or slowed down to ensure safety, thus generating further delays. If deviations are relatively small, they can be handled by modifying the timetable, with no changes to the duties for rolling stock and crew. On the other hand, massive disruptions of service (due, e.g., to strikes or technical issues) require both the timetable and the duties for rolling stock and crew to be modified. See \cite{cacchiani2014overview} for a review. Unlike the AGV applications our paper is motivated by, several minutes are usually allowed to generate a feasible solution which often makes feasible to use standard optimization software.  

Finally, our paper is related to the problem of \textit{coordinating the motions of multiple robots} operating in a shared workspace (e.g., robotic arms working in welding and painting workcells) without collisions. In this line of research, \cite{akella2002coordinating} study the problem when only the robot start times can be varied. They show that, even when the robot trajectories are specified, minimum time coordination of multiple robots is NP-hard, and define a Mixed Integer Linear Programming formulation. This work has been subsequently extended by \cite{peng2005coordinating} that determine velocity profiles for given paths that obey the kinematic constraints of robots
and avoid collisions while minimizing makespan. Unlike the fast exact algorithms presented in this paper, these authors use standard optimization software which make their approach impractical in large-scale real-time AGV applications.
For an up-to-date state-of-the-art of robot coordination, see \cite{spensieri2021modeling}.

As anticipated, in this paper we focus on recovering feasibility in \textit{Conflict-Free Vehicle Routing Problems} in industrial settings, like those arising in intralogistics applications \citep{fragapane2021planning}, where vehicles (AGVs) move at high speed (up to two meters per second) along relatively short ``segments" (e.g., aisles) of the guide path network. As a result, when deviations are observed, plans have to be revised in a few milliseconds to avoid collisions among vehicles. In this paper we identify two corrective actions and present fast exact algorithms to recover plan feasibility in real-time under four different objectives. In Section 2, we define the notation used throughout the paper and model the problem. In Section 3, we study the problem of recovering feasibility with proper vehicle corrective delays. In particular, we first present a model and then show that it can be cast as a shortest path problem on a suitably defined auxiliary graph for four distinct objective functions.
In Section 4, we extend the previous results to recovering feasibility with both vehicle delays and anticipations. In particular, we show that, under mild hypotheses, the problem can be solved with a two-stage approach in which a shortest path problem is solved at each stage.
Finally, in Section 5, we compare experimentally our tailored exact algorithms with solving the proposed formulations with off-the-shelf solvers.

\section{Problem definition and notation} \label{notation}
Let $G(V,A)$ be the graph representing the internal transportation network of a facility served by a fleet of vehicles. Vertex set $V$ represents loading/unloading stations, storage positions as well as intersections  while $A$ describes the segments of the material handling network. In practice, some complicating issues may be present to account for peculiar features of the vehicles (e.g., orientation as it is the case of automated forklifts). However, we neglect these aspects since they do not affect our findings.

In this article, we assume that a \textit{nominal} plan has been previously generated for a given \textit{Conflict-Free Vehicle Routing and Scheduling Problem} on $G$. As a result, each vehicle has been assigned a route and a schedule, including possible stops at some vertices in $V$. The individual vehicle plans constitute a nominal plan which is \textit{conflict-free}, and hopefully optimizes some performance measure. We also assume that the position of the vehicles is monitored periodically 
to identify possible delays or anticipations w.r.t. the nominal plan. If one or more vehicles come up to be early or late, the nominal plan may become unfeasible and some corrective actions need to be promptly identified and implemented in order to avoid collisions. 

In large-scale AGV applications, where computing time cannot exceed few milliseconds, feasibility recovery cannot include any route change and the only corrective actions are delaying some vehicles and possibly (depending on AGV technology) speeding up some vehicles.

Let $V_c=\{1, \dots, n\}$ be the set of vehicles and $\xi_h(t)$ the \textit{planned} position of vehicle $h \in V_c$ in $G$ at time $t$. Moreover, let $x_h(t)$ be the \textit{observed} position of vehicle $h \in V_c$ in $G$ at time $t$. Then, the deviation observed at time $t$ is 
$$d_h(t)=t-\xi_h^{-1}(x_h(t)).$$
It is worth noting that $d_h(t)$ can be either positive, negative or zero: if $d_h(t)>0$, vehicle $h$ is late w.r.t. the \textit{nominal} plan; if $d_h(t)<0$, it is ahead of its schedule; otherwise, it is on time. In what follows, for the sake of simplicity, we omit the time instant $t$ at which the AGV positions were last observed.

For every pair $h, k$ of vehicles, let $s_{h,k}$ ($\geq0$) be the maximum delay (\textit{slack}) that vehicle $h$ may accumulate, w.r.t. the nominal plan, without conflicting with vehicle $k$. Equivalently, $s_{h,k}$ represents the maximum conflict-free anticipation of vehicle $k$ on vehicle $h$ allowed by the current nominal plan. Slacks can be easily computed starting from the nominal plan. Finally, it is worth noting that slacks are not symmetric: e.g., $s_{h,k}=5$ and $s_{k,h}=+\infty$.

\begin{figure}[ht!]
\centering
\includegraphics[width=80mm]{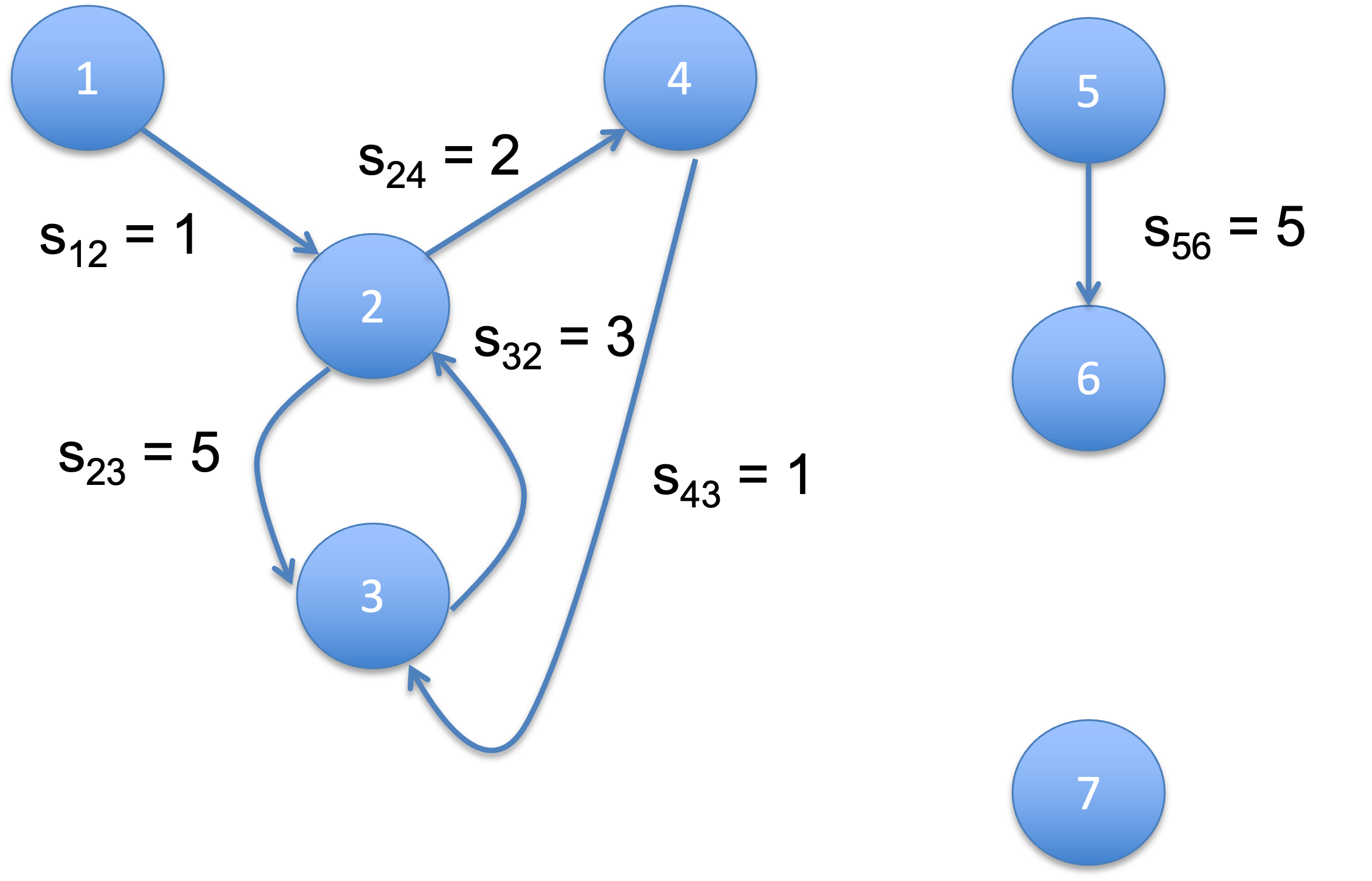}
\caption{A conflict graph \label{conflict_graph}}
\end{figure}

In order to model the relationship among vehicles in the nominal plan, we define a \textit{conflict graph} $G_c(V_c, A_c)$ where vertex $h \in V_c$ represents vehicle $h$ and an arc $(h,k)$ exists in $A_c$ if and only if $s_{h,k}<+\infty$. In this case the arc is assigned a weight equal to $s_{h,k}$. It is worth noting that graph $G_c$ may be cyclic and/or disconnected as shown in Figure \ref{conflict_graph}.

In order to avoid that deviations $d_h$ may cause vehicle conflicts, we impose to each vehicle $h \in V_c$ a delay or anticipation at the beginning of its route. Let $u_h  \in \mathbb{R}$ be the delay/anticipation of vehicle $h$ w.r.t. the nominal plan after feasibility has been recovered. Hence, 
$$\delta_h = u_h - d_h$$ 
represents the \textit{corrective} action imposed to vehicle $h$. If $\delta_h > 0$, a stop of duration $\delta_h$ will be inserted in the plan of the vehicle; if $\delta_h < 0$, vehicle $h$ will be accelerated in such a way as to generate an anticipation of $|\delta_h|$. When a vehicle $h$ does not respect its nominal plan ($d_h \neq 0$), it might be necessary to delay or anticipate other vehicles in order to recover feasibility. For instance, in Figure \ref{conflict_graph}, if $d_1=5$, $d_2=1$ and $d_3=\dots=d_7=0$, the initial delay of vehicles 1 and 2 may be compensated by further delaying vehicle 2 of 3 time units ($\delta_2=3$, $u_2=4$) which in turns requires vehicle 4 to be delayed of 2 time units ($\delta_4=u_4=2$) which ultimately requires vehicle 3 to be delayed of 1 time units ($\delta_3=u_3=1$). It is worth noting that vehicle 3 is delayed because of the interference between its route and that of vehicle 4 (not vehicle 2 because $s_{23}=5$ was greater than the difference between $u_2=4$ and $u_3=0$). Finally, vehicle 5, 6 and 7 were not affected by the initial plan perturbation.

\section{Recovering feasibility with corrective delays} \label{section_A}
We first model the problem of finding an optimal recovery plan under the hypothesis that only delays can be imposed to vehicles.  Assuming that the objective is to minimize the sum of vehicle delays, the problem can be formulated as follows:
\begin{align}
(P^D_1)\quad&\min z_1 = \sum\limits_{h \in V_c} u_h &&\label{1} \\
\mathrm{s.t.}\quad& &&\nonumber \\
& u_h - u_k \leq s_{h,k} & (h,k) \in A_c,& \label{2} \\
& u_h \geq d_h & h \in V_c,& \label{3}
\end{align}
where constraints (\ref{2}) impose that collisions between pairs of vehicles are avoided, while inequalities (\ref{3}) require that the total delay of a vehicle is no less that its observed deviation. 

A second objective is the minimization of a weighted sum of the vehicle delays. Let $w_h \geq 0$ be the cost of having vehicle $h$ delayed by a time unit. Then, the model can be reformulated as:
\begin{align}
(P^D_2)\quad&\min z_2 = \sum\limits_{h \in V_c} w_h u_h &&\nonumber \\
\mathrm{s.t. }\quad&  (\ref{2}) \textrm{ - } (\ref{3}).&&\nonumber 
\end{align}

A third objective is the minimization of the longest vehicle route (\textit{makespan}). Let $c_h$ be the completion time of vehicle $h$ in the nominal plan. Then, the model can be rewritten as:
\begin{align}
(P^D_3)\quad&\min z_3 = \max\limits_{h \in V_c} \left( c_h + u_h \right) &&\nonumber \\
\mathrm{s.t. }\quad& (\ref{2}) \textrm{ - } (\ref{3}),&&\nonumber 
\end{align}
or equivalently:
\begin{align}
(P^D_3)\quad&\min z_3  &&\nonumber \\
\mathrm{s.t. }\quad&  (\ref{2}) \textrm{ - } (\ref{3}), \textrm{ and }&&\nonumber \\
& z_3 \geq c_h + u_h  & h \in V_c.& \nonumber
\end{align}

Finally, we consider the case in which a \textit{due date} $\rho_h$ is specified for each vehicle $h \in V_c$ and the objective is the minimization of the total  \textit{lateness} of the vehicles w.r.t. their \textit{due dates}:
\begin{align}
(P^D_4)\quad&\min z_4  = \sum\limits_{h \in V_c} \max \left(0, u_h - \rho_h \right)&&\nonumber \\
\mathrm{s.t. }\quad& (\ref{2}) \textrm{ - } (\ref{3}),&&\nonumber 
\end{align}
or equivalently: 
\begin{align}
(P^D_4)\quad&\min z_4 = \sum\limits_{h \in V_c} y_h &&\nonumber \\
\mathrm{s.t. }\quad&  (\ref{2}) \textrm{ - } (\ref{3}), \textrm{ and }&&\nonumber \\
& y_h \geq u_h - \rho_h  & h \in V_c,&\nonumber\\
& y_h \geq 0 & h \in V_c.&\nonumber 
\end{align}

We now prove some propositions that will be helpful to derive tailored fast exact algorithms for recovering plan feasibility.
\begin{proposition}
Problems ($P^D_1$)-($P^D_4$) are feasible.
\end{proposition}
\begin{proof}
We prove this proposition for problem ($P^D_1$). Solution
\begin{eqnarray}
&& u_h = \max_{i \in V_c} d_i  \qquad h \in V_c, \label{4}
\end{eqnarray}
makes the left-hand sides of constraints (\ref{2}) null and satisfies constraints (\ref{3}). For problems ($P^D_2$)-($P^D_4$), the proof is similar.
\end{proof}

\noindent Solution (\ref{4}) is equivalent to impose to vehicle $h$ an additional delay equal to 
\begin{eqnarray}
&& \delta_h = \max_{i \in V_c} d_i - d_h. \label{5}
\end{eqnarray}

\noindent For instance, if the first vehicle is 1 unit of time late and the other vehicles are on time ($d_1=1$ and $d_h=0$ for $h=2,\dots, n)$, then a feasible solution (\ref{4}) is: $u_1 - d_1 = 0$ and $u_h - d_h = 1$ for $h=2,\dots, n$. On the other hand, if the first vehicle is one unit of time ahead and the other vehicles are on time ($d_1=-1$ and $d_h=0$ for $h=2,\dots, n$), then feasible solution (\ref{4}) is: $u_1 - d_1 = 1$ and $u_h - d_h=0$ for $h=2,\dots, n$.

\begin{proposition}\label{dualproblem}
Problem ($P^D_1$) is the dual of a one-to-all shortest path problem on a suitably defined auxiliary graph.
\end{proposition}
\begin{proof}
We define an auxiliary graph $\hat{G}_c = (\hat{V}_c,\hat{A}_c)$ with an additional (source) vertex $0$ and $|V_c|$ additional arcs from $0$ to any other node $h \in V_c$ having weight equal to $- d_h$. In other terms:
\begin{align*}
&\hat{V}_c =  V_c \cup \left\lbrace 0 \right\rbrace, & \\
&\hat{A}_c = A_c  \cup \left\lbrace (0, h) \ | \ h \in V_c \right\rbrace, & \\
&s_{0,h} = -d_h   & h \in V_c, \\
&s_{h,0} = +\infty   & h \in V_c. 
\end{align*}

Given a vertex $k \in \hat{V}_c$, we denote with $\delta^-(k)=\{h:(h,k) \in \hat{A}_c \}$ and $\delta^+(k)=\{h:(k,h) \in \hat{A}_c \}$ the set of in-neighbors and out-neighbors of $k$, respectively. 
Obviously: 
\begin{align*}
&|\delta^-(0)|=0,\\
&|\delta^+(0)|=|V_c|,\\
&u_0=d_0=0.
\end{align*}

Therefore, for any $k \in V_c$:
\begin{align}
    & u_{k} \geq u_{h} - s_{h,k} & h \in \delta^-(k),&\label{10}\\
    & u_{k} \geq d_k. &&\label{11}
\end{align}

Constraints (\ref{10}) and (\ref{11}) imply
\begin{align}
    & u_{k} \geq \max\limits_{h \in \delta^-(k)} \left( u_{h} - s_{h,k} \right) & k \in V_c. \label{recursive}
\end{align}

In order to minimize (\ref{1}), the recursive inequality (\ref{recursive}) must be satisfied as an equality, i.e.:
\begin{align}
    & u^*_{k} = \max\limits_{h \in \delta^-(k)} \left( u^*_{h} - s_{h,k} \right) & k \in V_c, \label{ref2}    
\end{align}
where the asterisk indicates an optimal solution.	
%It is worth noting that $|\delta^-(h)| \geq 1$ $h \in V_c$ (vertex $0$ is adjacent to any other vertex in $\hat{G}_c$).
For $h \neq 0$ relationship (\ref{ref2}) can be written as:
\begin{align}
    & u^*_{k} = \max\limits_{h \in \delta^-(k)} \left( \max\limits_{\ell \in \delta^-(h)} \left( u^*_{\ell} - s_{\ell, h} \right) - s_{h,k} \right) & k \in V_c.    \nonumber
\end{align}
By iterating until vertex $0$ is reached, the following expression is obtained:
\begin{align}\label{ref3}
    & u^*_{k} = \max\limits_{h \in V_c} \max\limits_{p \in \Pi_{hk}} \left( d_h - \sum\limits_{i = 1}^{|p|-1} s_{v_{i},v_{i+1}} \right) & k \in V_c,  \nonumber  
\end{align}
where $\Pi_{hk}$ ($(h,k \in V_c$) denotes the set of all paths in $G_c$ departing from $h \in V_c$ and arriving in $k \in V_c$, and $v_{i} \in V_c$ is the i-th vertex along path $p \in \Pi_{hk}$. Hence, the optimal solution $u_k^*$ corresponds to the shortest path between vertex $0$ and vertex $k$ on $\hat{G}_c$. 
\end{proof}

\begin{proposition}
Let $u_k^*$ ($k \in V_c$) be an optimal solution for problem ($P^D_1$). Then, $u_k^*$ ($k \in V_c$) is also optimal for problems ($P^D_2$)-($P^D_4$) (albeit obviously with different objective function values).
\end{proposition}
\begin{proof}
Let $u_k^*$ be a feasible solution minimizing $z_1$.

\noindent \textbf{Case \(z_2\)} - 
Because of (\ref{ref2}), solution $u_k^*$ remains optimal if weights $w_h$ are any nonnegative numbers.

\noindent \textbf{Case \(z_3\)} - 
We define a new graph $\bar{G}_c = (\bar{V}_c,\bar{A}_c)$ with an additional (sink) vertex $n+1$ as well as with $|V_c|$ additional arcs going from each vertex $h \in V_c$ to $n+1$ with weight equal to $\max\limits_{k \in V_c} c_k - c_h$:
\begin{align*}
&\bar{V}_c =  V_c \cup \left\lbrace n+1 \right\rbrace, & \\
&\bar{A}_c = A_c \cup \left\lbrace (h,n+1) \ | \ h \in V_c \right\rbrace, & \\
&s_{h,n+1} =  \max\limits_{k \in V_c} c_k - c_h & h \in V_c, \\
&s_{n+1,h} = +\infty   & h \in V_c. 
\end{align*}
By setting $d_{n+1}=0$, an optimal solution corresponds to the shortest path tree on $\bar{G}_c$ originating in $n+1$. It is worth noting that $u^*_{n+1}$ represents the minimum increase of the maximum completion time, i.e. $$\min z_3 = \max\limits_{k \in V_c} c_k + u^*_{n+1}.$$

\noindent \textbf{Case \(z_4\)} - A new graph $\tilde{G}_c = (\tilde{V}_c,\tilde{A}_c)$ is defined by adding (to $G_c$) $|V_c|$ additional (sink) vertices and the corresponding incoming arcs as follows:
\begin{align*}
&\tilde{V}_c = V_c \cup \left\lbrace n+1, \dots, 2 n \right\rbrace,& \\
&\tilde{A}_c = A_c \cup \left\lbrace (h,n+h) \ | \ h \in V_c \right\rbrace, & \\
&s_{h,n+h} =  \rho_h & h \in V_c, \\
&s_{n+h,h} =  +\infty & h \in V_c. 
\end{align*}
Then, by setting $d_{n+h}=0$ for $h \in V_c$, an optimal solution corresponds to a one-to-all shortest path tree on $\tilde{G}_c$. It is worth noting that $u^*_{n+h}$ represents the optimal value of variable $y_h$ in the definition of $z_4$, i.e.:
$$y^*_h= u^*_{n+h}.$$
\end{proof}

\section{Recovering feasibility with corrective delays and anticipations} \label{section_B}
In the previous section we examined how to recover the feasibility of a nominal plan by imposing corrective delays to vehicles.
In this section we consider the case that vehicles may also be accelerated, resulting in corrective anticipations with respect to the current plan. As before, the aim is to minimize four of the most common performance measures: (1) total vehicle delay, (2) total weighted delay, (3) maximum route duration (\textit{makespan}) and (4) total lateness w.r.t. specified \textit{due dates}.

For the sake of simplicity, we assume there are just two speed levels, a nominal speed $v_1$ and a higher speed 
\begin{equation}
v_2 = k \cdot v_1, \nonumber
\end{equation}
with $k > 1$, which can be hold for at most $T$ seconds (to limit energy consumption as well as to reduce the wear and tear of vehicles). 

Let $x_h$ be a nonnegative continuous variable representing the corrective anticipation imposed to vehicle $h$ to recover plan feasibility. Two factors impose an upper bound on $x_h$: $T$ and the time $t^c_h$ after which a conflict would occur if vehicle $h$ continued moving at speed $v_1$. 
As for the first factor, the maximum anticipation is achieved by using speed $v_2$ (instead of $v_1$) for $T$ instants and amounts to 
$$\dfrac{v_2 T}{v_1} - T = (k-1) T.$$
As for the second factor, the distance from the current position of vehicle $h$ to the first conflict on its route is $v_1 t^c_h$. So speed $v_2$ can be kept for at most 
$$\frac{v_1 t^c_h}{v_2}=\frac{t^c_h}{k}$$ 
instants, resulting into a maximum anticipation equal to:
$$t^c_h-\frac{t^c_h}{k}=\frac{k-1}{k}t^c_h.$$ 
Hence, the upper bound on $x_h$ is
$$x_h \leq L_h,$$ 
where
$$L_h = \min((k-1) T, \frac{k-1}{k}t^c_h).$$

The objective function is defined as the weighted sum, with parameters $\alpha > 0 $ and $\beta \geq 0$, of performance measure $z_i$ ($i= 1, \dots, 4$) and the total anticipation $\sum\limits_{h \in V_c} x_h$ imposed to vehicles. Problem $P^D_i$, ($i = 1, \dots, 4$) is then re-formulated as follows:
\begin{align}
(P^{AD}_i)\quad \min \quad & \begingroup
z'_i = \alpha z_i + \beta \sum\limits_{h \in V_c} x_h
\endgroup & \label{m3_start}\\
\mathrm{s.t. }\quad &&\nonumber\\
& u_{h} - x_{h} - u_{k} + x_{k} \leq s_{h,k} & (h, k) \in A_c, \label{m3_c1}\\
&u_h \geq d_h & h\in V_c,\label{m3_c2}\\
&0 \leq x_h \leq L_h & h\in V_c, \label{m3_end}
\end{align}
where constraints (\ref{m3_c1}) prevent collisions between pairs of vehicles, inequalities (\ref{m3_c2}) impose the observed delays/anticipations and constraints (\ref{m3_end}) define upper bounds on vehicle anticipations. We now prove some propositions.
\begin{proposition}
Problems $P^{AD}_i$ ($i = 1, \dots, 4$) are feasible.
\end{proposition}
\begin{proof}
Solution 
$x_h=0$ ($h \in V_c$) and $u_h = \max_{i \in V_c} d_i$ ($h \in V_c$)
satisfies all constraints.
\end{proof}

\begin{proposition}
Each basic feasible solution of formulation (\ref{m3_start})-(\ref{m3_end}) has either $\delta_h=u_h-d_h=0$ or $x_h=0$ (or both) for any $h \in V_c$.
\end{proposition}
\begin{proof}
Let $\mathbf{u}$ be the vector of $u_h$ associated to vehicles $h \in V_c$.
We prove the statement by contradiction. Let us assume that it exists an optimal solution with
$u_h - d_h \geq x_h > 0$ for a vehicle $h \in V_c$. 
A new feasible solution can be defined as follows:
\begin{align*}    
u^\prime_h=u_h-x_h,\\
x^\prime_h=0,
\end{align*}
Evaluating both solutions w.r.t. objective function $z'_i$ ($i = 1, \dots, n$), the following relationship is obtained:
\begin{align*}
&\alpha \cdot z_i(\mathbf{u}^\prime) + \beta \cdot 0 < \alpha \cdot z_i(\mathbf{u}) + \beta x_h & i=1, \dots, 4.
\end{align*}

Since the new feasible solution has a lower objective function value than the optimal one, the hypothesis is contradicted.
On the other hand, if an optimal solution is such that $x_h \geq u_h - d_h > 0$, then a new feasible solution can be defined as follows:
\begin{align*}    
u^\prime_h=d_h,\\
x^\prime_h=x_h-u_h+d_h.
\end{align*}
Once again, evaluating both solutions w.r.t. objective function $z'_i$ ($i = 1, \dots, n$),
\begin{align*}
&\alpha \cdot z_i(\mathbf{u}^\prime) + \beta \left( x_h-u_h+d_h \right) < \alpha \cdot z_i(\mathbf{u}) + \beta x_h & i=1, \dots, 4,
\end{align*}
contradicts the hypothesis.
\end{proof}

\begin{proposition}\label{profInverse}
If $L_h = +\infty$, problem $P^{AD}_i$ ($i = 1, \dots, 4$) with only anticipations is the dual of a one-to-all shortest path problem on a suitably defined graph.
\end{proposition}
\begin{proof}
For $u_h=d_h$, formulation (\ref{m3_start})-(\ref{m3_end}) becomes:
\begin{align}    
& \begingroup
\alpha z_i+\beta \sum\limits_{h \in V_c} d_{h} + \beta \min \sum\limits_{h \in V_c} \bar{x}_h
\endgroup &\label{m4_start} \\
\mathrm{s.t. }&&\nonumber\\
& \bar{x}_{k} - \bar{x}_{h} \leq s_{h,k} & (h, k)  \in A_c \label{m4_c1}\\
&\bar{x}_h \geq - d_h & h\in V_c\label{m4_end}
\end{align}
where $\bar{x}_h = x_h - d_h$. 
Therefore, formulation (\ref{m4_start})-(\ref{m4_end}) is equivalent to (\ref{1})-(\ref{3}) formulated on the reverse graph $G^{-1}_c = (V_c, A^{-1}_c)$ with $$A^{-1}_c = \{(k,h) \ | \ (h,k) \in A_c \},$$
and observed delays/anticipations $-d_h$. 
\end{proof}
The following proposition allows to decompose problem $P^{AD}_i$ ($i = 1, \dots, 4$) in case minimizing $z_i$ is more of a priority than minimizing the sum of the $x_h$ ($h \in V_c$) variables.
\begin{proposition}\label{2step}
If $\alpha \gg \beta \geq 0$ (in particular, if $\beta = 0$), problem $P^{AD}_i$ ($i = 1, \dots, 4$) decomposes into two independent (one-to-all shortest path) subproblems.
\end{proposition}
\begin{proof}
Initially, we prescribe that each vehicle is imposed an anticipation equal to the maximum allowed (e.g., $x_h=L_h$). Then formulation (\ref{m3_start})-(\ref{m3_end}) becomes:
\begin{align}    
& \begingroup
\alpha \min z_i + \beta \sum\limits_{h \in V_c} L_h
\endgroup &\label{m5_start} \\
\mathrm{s.t. }&&\nonumber\\
& \bar{u}_{h} - \bar{u}_{k} \leq s_{h,k} & (h, k)  \in A_c \label{m5_c1}\\
& \bar{u}_h \geq d_h - L_h & h\in V_c\label{m5_end}
\end{align}
where $\bar{u}_h=u_h-L_h$. Problem (\ref{m5_start})-(\ref{m5_end}) is equivalent to problem (\ref{1})-(\ref{3}) on a graph coincident with $G_c$, except that observed delays/anticipations are equal to $d_h - L_h$. Hence, it can be solved as a one-to-all shortest path problem. 
The corresponding optimal corrective delays $\delta_h=\bar{u}_h - d_h + L_h$ are feasible for (\ref{m3_start})-(\ref{m3_end}). Moreover, if $\beta = 0$ they are also optimal. Otherwise, this solution has to be corrected in order to minimize the total anticipation.
To do so, it has to be noted that, because of Proposition 5, the overall corrective delay determined by model (\ref{m5_start})-(\ref{m5_end}) is the optimal one also for (\ref{m3_start})-(\ref{m3_end}), i.e. $ \delta^*_h = \max \{0, \bar{u}_h-d_h \cancel{+L_h} \cancel{-L_h}\} $, while the overall corrective anticipation will be $x_h=\max \{0, \cancel{L_h}-\bar{u}_h+d_h \cancel{-L_h} \}$.
\\ 
However, because the anticipations $x_h$ ($h \in V_c$) prescribed in (\ref{m3_start})-(\ref{m3_end}) can be larger than the optimal ones, this formulation with $u_h=\delta^*_h + d_h$ becomes:
\begin{align}
&  \begingroup
\alpha z_i +\beta\sum\limits_{h \in V_c} (\delta^*_{h}+d_h) + \beta \min \sum\limits_{h \in V_c} \bar{x}_h
\endgroup && \label{m6_start}\\
\mathrm{s.t. }&&&\nonumber\\
& \bar{x}_{k} - \bar{x}_{h}  \leq s_{h,k} & (h, k)  \in A_c\label{m6_c1}\\
& \bar{x}_h \geq - \delta^*_h - d_h & h\in V_c\label{m6_end}
\end{align}
where corrective actions are $\bar{x}_h+\delta^*_h+d_h= x_h$.
Solution of model (\ref{m5_start})-(\ref{m5_end}) ensures feasibility of model (\ref{m6_start})-(\ref{m6_end}). Moreover, 
$$\bar{x}_h + \delta^*_h +d_h \leq L_h$$ is surely satisfied by that solution.
Then applying proposition \ref{profInverse} the thesis is proved.
\end{proof}
Summing up, if $\beta = 0$, only (\ref{m5_start})-(\ref{m5_end}) has to be solved. Otherwise, it is required to solve (\ref{m6_start})-(\ref{m6_end}) in a second stage. Both problems are equivalent to solving a one-to-all shortest path problem.

\section{Computational Results} \label{sec:results}
The aim of our computational experiments was to assess whether our approach can be valuable to recover plan feasibility in a real-time setting like those arising when using AGVs in manufacturing, transportation and logistics applications. In particular, we compared the performance of the tailored approach described in Sections (\ref{section_A}) and (\ref{section_B}) versus finding exact solutions to problems $P^D_1$-$P^D_4$ and $P^{AD}_1$-$P^{AD}_4$ with two off-the-shelf solvers: commercial solver IBM ILOG CPLEX 22.1 \citep{cplex22} and open-source solver SCIP 8.0 \citep{BestuzhevaEtal2021OO}.

All the experiments were run on a standalone Linux machine with an Intel Core i7 processor composed by $4$ cores clocked at $2.5$ GHz and equipped with $16$ GB of RAM. All algorithms have been coded in C++. Shortest paths on auxiliary graphs were determined by the Dijkstra's algorithm \citep{dijkstra1959note} with a Fibonacci heap min-priority queue.

We first describe the generation of input data, and then present the results.
Instances have been randomly generated with a number of vehicles $n = |V_c| = 50,100,150,200,250,300$. 
A parameter $p$, defined as  $$p=\dfrac{|A_c|}{n^2-n}-1,$$ was used to control the sparsity level of the conflict graph. 
In particular, $p$ was set equal to $\{0, 0.25, 0.50, 0.75\}$.
Vehicle (observed) deviations from the nominal plans $d_h$ ($n=1, \dots, n$) were randomly generated according to a uniform distribution in the [-10, 10] range. Vehicle slacks $s_{h,k}$ ($(h,k) \in A_c$) were uniformly generated in interval $[0,13]$. As far as $z_2$ is concerned, weights $w_h$ ($n=1, \dots, n$) were uniformly generated in $[0, 1]$. As for $z_3$, completion times $c_h$ ($n=1, \dots, n$) were derived from a uniform distribution in [100, 110]. Regarding objective function $z_4$, due dates $\rho_h$ ($n=1, \dots, n$) were uniformly generated in interval $[0,10]$. Finally, when corrective anticipations are allowed, we set $\alpha = 1000$ and $\beta = 1$.
For any pair $\left(n, p\right)$ we generated $10$ instances, for a total of $240$ instances. All test files are available at \url{https://tdrouting.com/cfpdp}.

The results of our experiments are reported in Tables \ref{tableResults_1} to \ref{tableResults_4} for objective functions $z_1$ to $z_4$. 
In all tables, the first two columns are self-explanatory and the values reported are averaged across all instances. The remaining column headings are as follows:
\begin{itemize}
\item $SP_D$: time (in milliseconds) spent by the tailored approach whenever only corrective delays are allowed;
\item $SP_{AD}$: time (in milliseconds) spent by the tailored approach whenever both corrective anticipations and corrective delays are allowed;
\item $CPLEX_D$: time (in milliseconds) spent by the IBM ILOG CPLEX 22.1 solver on problem $P^D_i$ ($i= 1, \dots, 4$));
\item $SCIP_D$: time (in milliseconds) spent by the SCIP 8.0 solver on problem $P^D_i$ ($i= 1, \dots, 4$));
\item $CPLEX_{AD}$: time (in milliseconds) spent by the IBM ILOG CPLEX 22.1 solver on problem  $P^{AD}_i$ ($i= 1, \dots, 4$));
\item $SCIP_{AD}$: time (in milliseconds) spent by the SCIP 8.0 solver on problem $P^{AD}_i$ ($i= 1, \dots, 4$));
\item $DEV_i$: percentage deviation between the $z_i$ values of the solutions with and without corrective anticipations allowed, i.e.
\begin{equation*}
DEV_i = 100 \cdot \dfrac{z^D_i - z^{AD}_i}{z^D_i},
\end{equation*}
where $z^{AD}_i$ and $z^D_i$ are the values of objective function $z_i$ when recovering feasibility with and without corrective anticipations, respectively.
\end{itemize}

Above all, computational results show that, in terms of computing time, our tailored exact approach outperformed both IBM ILOG CPLEX 22.1 and SCIP 8.0 solvers, independently of the objective function. Even for small fleet sizes ($n$=50), our procedures were able to recover feasibility in less than a millisecond while the CPLEX and SCIP solvers took something between 40 and 200 milliseconds, depending mainly on the sparsity of the conflict graph and less on the objective function. For moderate fleet sizes ($n$=100) the approach based on off-the-shelf solvers became no more viable (with a computing time between 100 and 600 milliseconds) while our approach was still able to recover feasibility in less than a millisecond. The performance gap became even more relevant for large fleets in which case the approach based on off-the-shelf solvers became impractical.

More in detail, Tables \ref{tableResults_1}-\ref{tableResults_4} show that, as a rule, computing times were lower for sparser conflict graphs ($p=0.75$). For instance, when minimizing $z_1$ with corrective delays, CPLEX took 612.07 milliseconds on average to obtain the optimal solution with $p=0.75$ versus 2487.02 milliseconds with $p=0$. This trend was confirmed for $z_2$, $z_3$ and $z_4$, with or without corrective anticipations.

Another aspect to be considered is that on average both $SP_D$ and $SP_{AD}$ were not affected by the objective function to be minimized. Moreover, $SP_{AD}$ was 2.5 times greater than $SP_D$, while still remaining very moderate (at most 0.48 millisecond for $n$=300).

From a managerial point of view, allowing vehicles to speed up to avoid conflicts was beneficial in terms of all objective functions. However, while the reductions in total delay, total weighted delay and lateness were significant (14.4 \%, 14.4 \%, 24.90 \% on average, respectively), the impact on makespan minimization was negligible (1.1 \%).

\begin{table}[H] 
\center
\captionsetup{font=normalsize}
\caption{Computational results - minimizing $z_1$ (total delay) \label{tableResults_1}}
\footnotesize
\begin{tabular}{|cc|ccc|ccc|c|}
\cline{3-8}
\multicolumn{2}{c|}{\textbf{}}     & \multicolumn{3}{c|}{\textbf{Corrective delays}} & \multicolumn{3}{c|}{\textbf{Corrective  anticipations and delays}} & \multicolumn{1}{c}{\textbf{}}\\
\hline
\multicolumn{2}{|c|}{\textbf{Instances}}     & \textbf{$SP_D$} & \textbf{$CPLEX_D$} & \textbf{$SCIP_D$}& \textbf{$SP_{AD}$} & \textbf{$CPLEX_{AD}$} & \textbf{$SCIP_{AD}$} &\textbf{DEV}\\
$p$            & $n$      & [ms] & [ms] & [ms] & [ms] & [ms] & [ms] & [\%] \\
\hline
0                           & 50         & 0.01     & 150.02      & 85.84             & 0.03       & 204.89         & 127.87   & 19.0  \\
                            & 100        & 0.02     & 577.64      & 306.15            & 0.08       & 606.52         & 489.50   & 16.5  \\
                            & 150        & 0.05     & 1316.73     & 679.56            & 0.13       & 1399.28        & 1476.53  & 13.8  \\
                            & 200        & 0.08     & 2347.16     & 1201.17           & 0.18       & 2713.41        & 3185.72  & 12.8  \\
                            & 250        & 0.12     & 3928.97     & 1901.89           & 0.26       & 4598.07        & 5718.10  & 11.0  \\
                            & 300        & 0.16     & 6601.63     & 2715.84           & 0.37       & 6702.86        & 8978.28  & 10.0  \\
\hline
0.25                          & 50         & 0.01     & 106.05      & 70.99             & 0.03       & 146.69         & 94.45    & 23.5  \\
                            & 100        & 0.03     & 443.90      & 231.14            & 0.08       & 477.93         & 404.10   & 14.8  \\
                            & 150        & 0.05     & 949.17      & 523.53            & 0.13       & 1049.86        & 1010.91  & 11.7  \\
                            & 200        & 0.07     & 1728.66     & 924.57            & 0.19       & 1875.75        & 1712.98  & 11.8  \\
                            & 250        & 0.10     & 2729.27     & 1452.90           & 0.27       & 3043.46        & 3495.32  & 12.9  \\
                            & 300        & 0.15     & 4058.68     & 2162.66           & 0.39       & 4264.78        & 4468.80  & 9.9  \\
\hline                    
0.50                          & 50         & 0.01     & 90.15       & 53.89             & 0.03       & 84.00          & 67.71    & 26.4  \\
                            & 100        & 0.03     & 324.10      & 213.07            & 0.09       & 318.45         & 254.37   & 16.0  \\
                            & 150        & 0.06     & 683.22      & 374.80            & 0.14       & 722.09         & 765.84   & 10.9  \\
                            & 200        & 0.09     & 1193.58     & 647.02            & 0.21       & 1297.74        & 1282.30  & 13.0  \\
                            & 250        & 0.13     & 2056.36     & 1111.71           & 0.31       & 1937.82        & 1888.44  & 9.9  \\
                            & 300        & 0.16     & 2935.56     & 1595.95           & 0.43       & 2955.88        & 2785.92  & 11.8  \\
\hline
0.75                          & 50         & 0.01     & 45.33       & 30.39             & 0.03       & 51.35          & 42.16    & 25.8  \\
                            & 100        & 0.04     & 150.79      & 106.25            & 0.10       & 188.10         & 146.78   & 15.4  \\
                            & 150        & 0.07     & 342.28      & 230.57            & 0.18       & 405.06         & 358.31   & 15.8  \\
                            & 200        & 0.10     & 622.10      & 431.28            & 0.25       & 729.30         & 813.81   & 10.1  \\
                            & 250        & 0.17     & 1071.23     & 678.36            & 0.38       & 1139.11        & 1299.99  & 10.0  \\
                            & 300        & 0.19     & 1440.69     & 945.51            & 0.48       & 1556.97        & 1921.71  & 12.1  \\
\hline
\multicolumn{2}{|c|}{\textbf{AVERAGE}} & \textbf{0.08}  & \textbf{1495.55} & \textbf{778.13}   & \textbf{0.20}   & \textbf{1602.88}  & \textbf{1782.91}     & \textbf{14.4}             \\
\hline
\end{tabular}
\end{table}

\begin{table}[H] 
\center
\captionsetup{font=normalsize}
\caption{Computational results - minimizing $z_2$ (total weighted delay) \label{tableResults_2}}
\footnotesize
\begin{tabular}{|cc|ccc|ccc|c|}
\cline{3-8}
\multicolumn{2}{c|}{\textbf{}}     & \multicolumn{3}{c|}{\textbf{Corrective delays}} & \multicolumn{3}{c|}{\textbf{Corrective  anticipations and delays}} & \multicolumn{1}{c}{\textbf{}}\\
\hline
\multicolumn{2}{|c|}{\textbf{Instances}}     & \textbf{$SP_D$} & \textbf{$CPLEX_D$} & \textbf{$SCIP_D$}& \textbf{$SP_{AD}$} & \textbf{$CPLEX_{AD}$} & \textbf{$SCIP_{AD}$} &\textbf{DEV}\\
$p$            & $n$      & [ms] & [ms] & [ms] & [ms] & [ms] & [ms] & [\%] \\
\hline
0                    & 50       & 0.01                                             & 140.65                                             & 90.40                                              & 0.03                                    & 151.53                                             & 121.86                                             & 19.1                                               \\
                     & 100      & 0.02                                             & 570.53                                             & 295.62                                             & 0.08                                    & 594.22                                             & 469.01                                             & 16.0                                               \\
                     & 150      & 0.05                                             & 1306.04                                            & 656.59                                             & 0.13                                    & 1359.35                                            & 1400.32                                            & 14.0                                               \\
                     & 200      & 0.08                                             & 2367.73                                            & 1175.33                                            & 0.18                                    & 2455.25                                            & 2536.37                                            & 13.0                                               \\
                     & 250      & 0.12                                             & 3965.93                                            & 1853.67                                            & 0.26                                    & 3901.06                                            & 4474.88                                            & 11.0                                               \\
                     & 300      & 0.16                                             & 6656.24                                            & 2673.13                                            & 0.37                                    & 5685.01                                            & 7172.66                                            & 10.0                                               \\ \hline
0.25                   & 50       & 0.01                                             & 115.33                                             & 70.60                                              & 0.03                                    & 119.47                                             & 87.59                                              & 23.7                                               \\
                     & 100      & 0.03                                             & 439.93                                             & 232.99                                             & 0.08                                    & 451.78                                             & 375.09                                             & 15.0                                               \\
                     & 150      & 0.05                                             & 1015.91                                            & 507.37                                             & 0.13                                    & 1047.64                                            & 938.06                                             & 12.0                                               \\
                     & 200      & 0.07                                             & 1769.17                                            & 910.15                                             & 0.19                                    & 1836.91                                            & 1597.90                                            & 12.0                                               \\
                     & 250      & 0.10                                             & 2787.87                                            & 1406.95                                            & 0.27                                    & 2892.88                                            & 3031.46                                            & 13.0                                               \\
                     & 300      & 0.15                                             & 4024.27                                            & 2038.38                                            & 0.39                                    & 4159.17                                            & 4226.66                                            & 10.0                                               \\ \hline
0.50                   & 50       & 0.01                                             & 108.66                                             & 52.40                                              & 0.03                                    & 112.88                                             & 63.06                                              & 26.3                                               \\
                     & 100      & 0.03                                             & 289.31                                             & 180.71                                             & 0.09                                    & 302.22                                             & 248.92                                             & 15.1                                               \\
                     & 150      & 0.06                                             & 652.11                                             & 361.32                                             & 0.14                                    & 684.88                                             & 697.68                                             & 11.0                                               \\
                     & 200      & 0.09                                             & 1169.77                                            & 616.04                                             & 0.21                                    & 1228.86                                            & 1164.55                                            & 13.0                                               \\
                     & 250      & 0.13                                             & 1844.79                                            & 960.09                                             & 0.31                                    & 1938.11                                            & 1830.03                                            & 10.0                                               \\
                     & 300      & 0.16                                             & 2670.26                                            & 1382.90                                            & 0.43                                    & 2798.06                                            & 2447.01                                            & 12.0                                               \\ \hline
0.75                   & 50       & 0.01                                             & 40.88                                              & 27.18                                              & 0.03                                    & 43.21                                              & 36.93                                              & 25.6                                               \\
                     & 100      & 0.04                                             & 139.51                                             & 100.16                                             & 0.10                                    & 146.83                                             & 134.53                                             & 15.9                                               \\
                     & 150      & 0.07                                             & 332.19                                             & 224.70                                             & 0.18                                    & 349.67                                             & 312.36                                             & 16.1                                               \\
                     & 200      & 0.10                                             & 587.61                                             & 409.05                                             & 0.25                                    & 616.49                                             & 646.27                                             & 10.0                                               \\
                     & 250      & 0.14                                             & 915.13                                             & 580.78                                             & 0.38                                    & 961.58                                             & 1133.95                                            & 10.0                                               \\
                     & 300      & 0.19                                             & 1351.54                                            & 818.13                                             & 0.48                                    & 1394.44                                            & 1538.65                                            & 12.0                                               \\ \hline
\multicolumn{2}{|c|}{AVERAGE}       & \textbf{0.08}                                    & \textbf{1469.22}                                   & \textbf{734.36}                                    & \textbf{0.20}                           & \textbf{1467.98}                                   & \textbf{1528.57}                                   & \textbf{14.4}                                      \\ \hline
\end{tabular}
\end{table}

\begin{table}[H] 
\center
\captionsetup{font=normalsize}
\caption{Computational results - minimizing $z_3$ (maximum route duration) \label{tableResults_3}}
\footnotesize
\begin{tabular}{|cc|ccc|ccc|c|}
\cline{3-8}
\multicolumn{2}{c|}{\textbf{}}     & \multicolumn{3}{c|}{\textbf{Corrective delays}} & \multicolumn{3}{c|}{\textbf{Corrective  anticipations and delays}} & \multicolumn{1}{c}{\textbf{}}\\
\hline
\multicolumn{2}{|c|}{\textbf{Instances}}     & \textbf{$SP_D$} & \textbf{$CPLEX_D$} & \textbf{$SCIP_D$}& \textbf{$SP_{AD}$} & \textbf{$CPLEX_{AD}$} & \textbf{$SCIP_{AD}$} &\textbf{DEV}\\
$p$            & $n$      & [ms] & [ms] & [ms] & [ms] & [ms] & [ms] & [\%] \\
\hline
0                    & 50       & 0.01                                           & 141.80                                             & 87.19                                              & 0.03                                     & 157.64                                             & 92.15                                              & 1.5                                                \\
                     & 100      & 0.02                                           & 582.44                                             & 297.14                                             & 0.08                                     & 604.55                                             & 345.81                                             & 1.3                                                \\
                     & 150      & 0.05                                           & 1305.45                                            & 670.45                                             & 0.13                                     & 1356.92                                            & 954.05                                             & 1.2                                                \\
                     & 200      & 0.08                                           & 2357.89                                            & 1184.54                                            & 0.18                                     & 2444.40                                            & 1655.04                                            & 1.1                                                \\
                     & 250      & 0.12                                           & 3977.72                                            & 1862.65                                            & 0.26                                     & 3848.72                                            & 2558.27                                            & 0.9                                                \\
                     & 300      & 0.16                                           & 6697.59                                            & 2692.51                                            & 0.37                                     & 5495.73                                            & 3679.66                                            & 0.8                                                \\ \hline
0.25                   & 50       & 0.01                                           & 134.99                                             & 72.70                                              & 0.03                                     & 141.89                                             & 82.47                                              & 1.9                                                \\
                     & 100      & 0.03                                           & 440.30                                             & 238.24                                             & 0.08                                     & 446.43                                             & 316.78                                             & 1.3                                                \\
                     & 150      & 0.05                                           & 967.83                                             & 511.27                                             & 0.13                                     & 1006.07                                            & 895.67                                             & 1.0                                                \\
                     & 200      & 0.07                                           & 1769.53                                            & 914.28                                             & 0.19                                     & 1834.56                                            & 1488.75                                            & 1.0                                                \\
                     & 250      & 0.10                                           & 2776.47                                            & 1421.74                                            & 0.27                                     & 2880.65                                            & 2420.99                                            & 1.1                                                \\
                     & 300      & 0.15                                           & 4014.25                                            & 2066.91                                            & 0.39                                     & 4173.16                                            & 3453.27                                            & 0.8                                                \\ \hline
0.50                   & 50       & 0.01                                           & 74.39                                              & 54.60                                              & 0.03                                     & 77.33                                              & 64.33                                              & 1.9                                                \\
                     & 100      & 0.03                                           & 298.07                                             & 184.66                                             & 0.09                                     & 308.33                                             & 245.71                                             & 1.2                                                \\
                     & 150      & 0.06                                           & 659.08                                             & 358.32                                             & 0.14                                     & 687.23                                             & 659.39                                             & 0.9                                                \\
                     & 200      & 0.09                                           & 1180.71                                            & 620.85                                             & 0.21                                     & 1226.90                                            & 1110.35                                            & 1.1                                                \\
                     & 250      & 0.13                                           & 1866.16                                            & 967.59                                             & 0.31                                     & 1939.64                                            & 1695.69                                            & 0.8                                                \\
                     & 300      & 0.16                                           & 2672.27                                            & 1391.01                                            & 0.43                                     & 2792.27                                            & 2502.10                                            & 1.0                                                \\ \hline
0.75                   & 50       & 0.01                                           & 61.83                                              & 30.04                                              & 0.03                                     & 63.28                                              & 37.07                                              & 1.2                                                \\
                     & 100      & 0.04                                           & 196.10                                             & 105.34                                             & 0.10                                     & 205.19                                             & 138.50                                             & 1.3                                                \\
                     & 150      & 0.07                                           & 342.35                                             & 237.18                                             & 0.18                                     & 352.69                                             & 321.32                                             & 1.2                                                \\
                     & 200      & 0.10                                           & 596.75                                             & 398.47                                             & 0.25                                     & 617.89                                             & 655.02                                             & 0.8                                                \\
                     & 250      & 0.14                                           & 928.92                                             & 608.53                                             & 0.38                                     & 959.61                                             & 1052.39                                            & 0.8                                                \\
                     & 300      & 0.19                                           & 1366.83                                            & 842.93                                             & 0.48                                     & 1387.70                                            & 1650.93                                            & 1.0                                                \\ \hline
\multicolumn{2}{|c|}{AVERAGE}       & \textbf{0.08}                                  & \textbf{1475.40}                                   & \textbf{742.46}                                    & \textbf{0.20}                            & \textbf{1458.70}                                   & \textbf{1169.82}                                   & \textbf{1.1}                                       \\ \hline
\end{tabular}
\end{table}

\begin{table}[H] 
\center
\captionsetup{font=normalsize}
\caption{Computational results - minimizing $z_4$ (total lateness) \label{tableResults_4}}
\footnotesize
\begin{tabular}{|cc|ccc|ccc|c|}
\cline{3-8}
\multicolumn{2}{c|}{\textbf{}}     & \multicolumn{3}{c|}{\textbf{Corrective delays}} & \multicolumn{3}{c|}{\textbf{Corrective  anticipations and delays}} & \multicolumn{1}{c}{\textbf{}}\\
\hline
\multicolumn{2}{|c|}{\textbf{Instances}}     & \textbf{$SP_D$} & \textbf{$CPLEX_D$} & \textbf{$SCIP_D$}& \textbf{$SP_{AD}$} & \textbf{$CPLEX_{AD}$} & \textbf{$SCIP_{AD}$} &\textbf{DEV}\\
$p$            & $n$      & [ms] & [ms] & [ms] & [ms] & [ms] & [ms] & [\%] \\
\hline
0                    & 50       & 0.01                                             & 153.23                                             & 87.77                                              & 0.03                                             & 163.045                                             & 95.94                                              & 31.5                                               \\
                     & 100      & 0.02                                             & 588.34                                             & 300.91                                             & 0.08                                             & 617.637                                             & 347.33                                             & 26.5                                               \\
                     & 150      & 0.05                                             & 1319.74                                            & 668.03                                             & 0.13                                             & 1380.705                                            & 947.73                                             & 25.0                                               \\
                     & 200      & 0.08                                             & 2373.19                                            & 1192.34                                            & 0.18                                             & 2467.960                                            & 1663.02                                            & 23.1                                               \\
                     & 250      & 0.12                                             & 3984.24                                            & 1871.59                                            & 0.26                                             & 3871.128                                            & 2572.47                                            & 20.3                                               \\
                     & 300      & 0.16                                             & 6749.59                                            & 2707.37                                            & 0.37                                             & 5574.294                                            & 3674.64                                            & 18.3                                               \\ \hline
0.25                   & 50       & 0.01                                             & 113.05                                             & 71.92                                              & 0.03                                             & 117.086                                             & 87.87                                              & 39.2                                               \\
                     & 100      & 0.03                                             & 458.40                                             & 237.54                                             & 0.08                                             & 466.384                                             & 326.79                                             & 25.9                                               \\
                     & 150      & 0.05                                             & 1020.14                                            & 516.81                                             & 0.13                                             & 1056.351                                            & 908.94                                             & 20.9                                               \\
                     & 200      & 0.07                                             & 1757.49                                            & 944.92                                             & 0.19                                             & 1831.170                                            & 1682.39                                            & 21.6                                               \\
                     & 250      & 0.10                                             & 2785.57                                            & 1434.06                                            & 0.27                                             & 2913.645                                            & 3283.40                                            & 23.2                                               \\
                     & 300      & 0.15                                             & 4020.14                                            & 2070.30                                            & 0.39                                             & 4194.527                                            & 4055.17                                            & 18.3                                               \\ \hline
0.50                   & 50       & 0.01                                             & 74.26                                              & 54.09                                              & 0.03                                             & 78.165                                              & 68.48                                              & 44.3                                               \\
                     & 100      & 0.03                                             & 295.42                                             & 184.98                                             & 0.09                                             & 306.283                                             & 266.42                                             & 26.2                                               \\
                     & 150      & 0.06                                             & 660.25                                             & 363.06                                             & 0.14                                             & 691.131                                             & 708.29                                             & 20.7                                               \\
                     & 200      & 0.09                                             & 1190.09                                            & 626.48                                             & 0.21                                             & 1246.885                                            & 1186.81                                            & 23.2                                               \\
                     & 250      & 0.13                                             & 1849.15                                            & 980.66                                             & 0.31                                             & 1939.435                                            & 1908.06                                            & 18.2                                               \\
                     & 300      & 0.16                                             & 2696.32                                            & 1405.35                                            & 0.43                                             & 2829.059                                            & 2654.74                                            & 21.5                                               \\ \hline
0.75                   & 50       & 0.01                                             & 43.19                                              & 29.39                                              & 0.03                                             & 45.024                                              & 43.01                                              & 38.0                                               \\
                     & 100      & 0.04                                             & 143.70                                             & 102.76                                             & 0.10                                             & 151.063                                             & 151.51                                             & 27.1                                               \\
                     & 150      & 0.07                                             & 342.50                                             & 232.38                                             & 0.18                                             & 356.111                                             & 345.88                                             & 27.3                                               \\
                     & 200      & 0.10                                             & 600.56                                             & 414.04                                             & 0.25                                             & 626.077                                             & 712.76                                             & 17.8                                               \\
                     & 250      & 0.14                                             & 934.83                                             & 602.25                                             & 0.38                                             & 969.996                                             & 1143.79                                            & 18.4                                               \\
                     & 300      & 0.19                                             & 1367.17                                            & 835.19                                             & 0.48                                             & 1401.898                                            & 1708.89                                            & 21.3                                               \\ \hline
\multicolumn{2}{|c|}{AVERAGE}       & \textbf{0.08}                                    & \textbf{1480.02}                                   & \textbf{747.26}                                    & \textbf{0.20}                                    & \textbf{1470.63}                                   & \textbf{1272.68}                                   & \textbf{24.90}                                      \\ \hline
\end{tabular}
\end{table}

\section{Conclusions}
\label{sec:conclusions}
In this paper we have dealt with recovering plan feasibility in \textit{Conflict-Free Vehicle Routing and Scheduling Problems} whenever some vehicles are ahead or behind of schedule. The problem is of the outmost importance in manufacturing, transportation and logistics facilities where AGVs are utilized to move loads between stations. 
In such settings, vehicles move at speeds in the order of one-two meter per second and feasibility has to be recovered in a few milliseconds. In this paper we have presented fast exact algorithms to solve this problem with respect to two corrective actions (introducing corrective delays/anticipations) with the objective to optimize four common performance measures (total vehicle delay, total weighted delay, maximum route duration and total lateness). An extensive empirical study has shown that, in terms of computing time, our tailored exact algorithms are at least three orders of magnitude faster than IBM ILOG CPLEX 22.1 \citep{cplex22} and SCIP 8.0 \citep{BestuzhevaEtal2021OO} solvers and are suitable for large intralogistics applications.

 \section*{Conflict of interest}
 The authors declare that they have no conflict of interest.


%\bibliographystyle{model5-names}\biboptions{authoryear}
\bibliography{biblio.bib}
%% Non-BibTeX users please use
%\begin{thebibliography}{}
%%
%% and use \bibitem to create references. Consult the Instructions
%% for authors for reference list style.
%\end{thebibliography}

\end{document}